# On Evaluation of Nonlinear Exponential Sums
N. A. Carella, June 2002

**Abstract:** This paper provides a technique for evaluating the nonlinear exponential sums $T_d(\chi,s) = \Sigma_{0 \leq x < p}\chi(x)\exp(i2\pi s x^d/p)$ in closed forms. The evaluation is obtained from the known values of simpler exponential sums.

## 1 Introduction

This paper provides a technique for evaluating the nonlinear exponential sums $T_d(\chi,s) = \Sigma_{0 \leq x < p}\chi(x)\exp(i2\pi s x^d/p)$ in closed forms. The evaluations are obtained from the known values of simpler exponential sums. Investigation of the properties of the related gaussian sums $S_k(s) = \Sigma_{0 \leq x < p}\exp(i2\pi s x^k/p)$ and $G_n(\chi,s) = \Sigma_{0 \leq x < p}\chi(x)\exp(i2\pi s x/p)$ has a long history, it extends back to at least two centuries. Very elaborate algebraic and analytical techniques have been developed to deal with this problem. Nonetheless, many questions about these sums are still unresolved, even for the simplest cases of the cubic, quartic and quintic gaussian sums.

## 2 Closed Forms Evaluations of $T_d(\chi,s)$

Let $\chi$ be a multiplicative character on $\mathbf{F}_p$ of order $\mathrm{ord}(\chi) = n$. The exact value of the nonlinear exponential sum

$$(1) \qquad T_d(\chi,s) = \sum_{x=0}^{p-1} \chi(x) e^{i2\pi s x^d / p}$$

of degree $d$ will be determined in terms of the values of sums $S_2, S_3, \ldots, S_k$, or equivalently as a polynomial $T_d(\chi,s) \in \mathbb{Z}[S_2, S_3, \ldots, S_k]$. The evaluations already known are the followings:

(1) $T_1(\chi,s) = G_n(\chi,s)$, with characters of order $n = \mathrm{ord}(\chi) = 2, 3$, and 4. The analysis of the cases $n = 3, 4$ appear in [6], [7], see also [3], [5], [9], and [10], etc.
(2) $T_d(\chi,s) = 0$, with $\mathrm{ord}(\chi) = 2$, $d$ even and $p \equiv 3 \bmod 4$.

Case (2) above is easy. The remaining sums $T_d(\chi,s)$ with $p - 1 \equiv 0 \bmod d$ are not necessarily trivial. The closed form evaluation of $T_2(\chi)$, $T_3(\chi)$, and $T_4(\chi)$, $\text{ord}(\chi) = 2$, up to a sign will be determined next. Some of the sign ambiguity has not been resolved yet.

**Theorem 1.** Let $p = a_4^2 + b_4^2$ be a prime for which $a_4 \equiv -\mu \bmod 4$, $\mu \equiv 2^{(p-1)/2} \bmod p$, and let $\text{ord}(\chi) = 2$. Then

(2) $$\sum_{x=0}^{p-1} \left(\frac{x}{p}\right) e^{i2\pi x^2/p} = \pm\sqrt{2\mu(p + a_4\sqrt{p})}.$$

Proof: It is simple matter to show that $S_4 = S_2 + T_2(\chi)$. Substituting the value of $S_4$, [2, p. 362], completes the claim. (Note that $T_2(\chi) = 0$ for $p = 4m + 3$). ∎

The symbol $i^* = (-1)^{(p-1)^2/8}$ used below is equal to either 1 or $i$ depending on whether $p \equiv 1$ or 3 mod 4.

**Theorem 2.** Let $p \equiv 1 \bmod 6$, and let $\text{ord}(\chi) = 2$. Then

(3) $$T_3(\chi) = \begin{cases} i^*(S_3^2 - p)p^{-1/2} & \text{if 2 is a cubic residue,} \\ i^*(4p - S_3^2 \pm S_3(12p - 3S_3^2)^{1/2})p^{-1/2}/2 & \text{if 2 is a cubic nonresidue.} \end{cases}$$

Proof: Consider the equation $S_6 = S_3 + T_3(\chi)$, and substitute the value of $S_6$, [2, p. 360]. ∎

The case of $p \equiv 5 \bmod 6$ reduces to $S_6 = S_2$, this follows from $S_k = S_d$, for $d = \gcd(k, p - 1)$. A theorem of Gauss asserts that $2 \in \mathbf{F}_p$ is a cubic residue if and only if the prime $p = 3n + 1$ has the quadratic form representation $p = x^2 + 27y^2$, for some integers $x, y \in \mathbb{Z}$. So the complex number $T_3(\chi)$ assumes the first entry in the formula if $p$ is a prime of this form. The sign in the second part of the formula is not difficult to determine, for details, see [2, p. 360].

**Corollary 3.** If $\chi$ is the quadratic symbol and $p \equiv 1 \bmod 6$, then

(4i) $\sqrt{p} \leq |T_3(\chi,s)| \leq 3\sqrt{p} \quad \Leftrightarrow \quad p = x^2 + 3y^2,$

(4ii) $0 \leq |T_3(\chi,s)| \leq 3\sqrt{p} \quad \Leftrightarrow \quad p \neq x^2 + 3y^2.$

**Theorem 4.** Let $p = a_8^2 + 2b_8^2 = 8m + 1$ be prime, $a_8 \equiv -1 \bmod 4$, and $\eta \equiv 2^{(p-1)/4}$ mod $p$. Then

$$\text{(5)} \quad \sum_{x=0}^{p-1} \chi(x) e^{i2\pi x^4/p} = \pm\sqrt{(a_8 + \sqrt{p})(2\eta\sqrt{2(p + a_4\sqrt{p})} + (-1)^m \sqrt{p})}.$$

The proof/evaluation of the sum $T_d(\chi)$, ord($\chi$) = 2, has the same format: $T_d(\chi,s) = S_{2d}(s) - S_d(s)$, with gcd$(2d, p-1) = 2d$. Currently it can also be evaluated for $d = 2, 3, 4, 6, 8$ and $12$ in terms of the known values $S_2, S_4, S_6, S_8, S_{12}, S_{16}$, and $S_{24}$, up to signs see [2], [1]. In addition, the exact values of $S_k$ over $\mathbf{F}_q$ for the parameter $k = p^e + 1$, $p \mid q$, are given in [4].

The last sum considered is evaluated in terms of the well known Kloostermann sum and Salie sum:

$$\text{(6)} \quad K(a,b) = \sum_{x=0}^{p-1} e^{i2\pi(ax + bx^{-1})/p} \quad \text{and} \quad S(a,b) = \sum_{x=0}^{p-1} \left(\frac{x}{p}\right) e^{i2\pi(ax + bx^{-1})/p},$$

respectively

**Theorem 5.** Let $p$ be a prime, and let $0 \neq a, b \in \mathbb{Z}$. Then

$$\text{(7)} \quad \sum_{x=0}^{p-1} e^{i2\pi(ax^2 + bx^{-2})/p} = 2\sqrt{p}\left(\left(\frac{a}{p}\right)(-1)^{(p-1)^2/8} \cos(4\pi a/p) + \cos(2\pi\theta/p)\right),$$

where the $p - 1$ discrete angles $\theta \in [0, \pi]$ are distributed in accordance to the Sato-Tate measure $2\pi^{-1}\sin(\theta)d\theta$.